# Regional Boundary Asymptotic Gradient Full Order Observer Via Internal Region

Raheam A. Al-Saphory[1] and Zinah A. Khalid[2]

[1*, 2] Department of Mathematics, College of Education for Pure Sciences, Tikrit University, Tikrit, Iraq.

Email:(saphory@hotmail.com;zeeena_assif@hotmail.com)

*Corresponding Author

**Abstract:**

The purpose of this paper is to deals with the problem of regional boundary asymptotic gradient full order observer ($\mathbf{\Gamma^*AGFO}$-observer) concept by using internal regional case. Thus, we study the relation between this notion and the corresponding asymptotic detectability and sensors. More precisely, various important results have been examined and explored concern an extension of an approach which enables to reconstruct the gradient of current state from internal region. In addition, it has been shown that the characterization of $\mathbf{\Gamma^*AGFO}$-observability under which conditions to be achieved. Finally, we have illustrate that there is a dynamical system which does not represent the observer in the usual sense, but it could be interpreted as a $\mathbf{\Gamma^*AGFO}$-observer.

## 1. Introduction

The concept of asymptotic observer theory was discovered by Luenberger for finite dimensional linear system as in [1]. Thus, this approach has been generalized to distributed parameter systems characterized by strongly continuous semi-group operators in Hilbert space by Grassing and Lamont [2]. The characterization of an asymptotic observer via sensor and actuator structures was explored by El Jai *et al.* in [3-5].

Recently, an important extension is that of regional and regional boundary state reconstruction has been introduced by Zerrik and El Jai *et al.* for finite time [6-9]. Al-saphory and El-Jai introduce the asymptotic regional state observation in infinite time as in [10-12].Therefore, the regional analysis consists in studying the asymptotic behavior of the systems not in the whole the domain but only in region $\omega \subset \Omega$ or on $\Gamma \subset \partial\Omega$ of system domain $\Omega$ [13-16].

Another orientation of regional gradient analysis in more different systems and regions [17-19] and for asymptotic case [20-22].

The purpose of this concept is motivated by certain concrete-real problems, in mechanic, thermic, environment in [23-25]. In this paper, we explore an approach which allow to construct $\Gamma^*$AGFO-observer in a given region $\Gamma^*$of the domain boundary $\partial\Omega$in connection with regional boundary gradient strategic sensor ($\Gamma^*$G-strategic sensor)and regional boundary asymptotic gradient detectability ($\Gamma^*$AG-detectability).

It is interested to study the problem of the treatment of water by using a bioreactor where the objective is to estimate the concentration of substrate at the boundary output of the bioreactor in order the water regulation is achieved ( figure 1) [26].

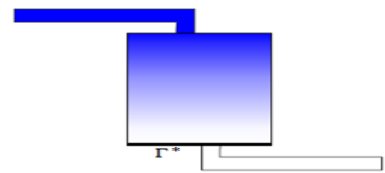

**Fig. 1:** Substrate concentration at the output of the reactor $\Gamma^*$.

The outline of this paper is organized as follow:

Section 2 concerns the class of considered system, definition, characterizations in connection with sensors and preliminaries of regional boundary gradient observability and detectability. Section 3, devotes to the problem of crossing method from internal region to boundary case by using trace operator esteems. Section 4, gives an application to various situations of sensors locations on the regional boundary gradient detectability in diffusion parabolic distributed systems. Last section tackles the relation between regional boundary detectability of state gradient and regional boundary observer.

## 2. Considered System and Problem Formulation

Consider a distributed parameter system defined with the following forms:
⋄ $\Omega$ is an open bounded subset of $R^n$with smooth boundary $\partial\Omega$.
⋄ $\Gamma$ is a sub-region of $\partial\Omega$ with positive measure.
⋄ Denote $Q = \Omega \times ]0,\infty[$ and $\Theta = \partial\Omega \times ]0,\infty[$.
⋄ The space $X = H^1(\Omega)$, $U = L^2(0,T,R^p)$ and $\mathcal{O} = L^2(0,T,R^q)$ are designed in this paper as separable Hilbert spaces and represented as state space, control space and observation space where $p$ and $q$ are the numbers of actuators and sensors [10].

◇ $A = \sum_{i,j=1}^{n} \frac{\partial}{\partial x_j}(a_{ij}\frac{\partial}{\partial x_j})$ with $a_{ij} \in D(\bar{A})$ (domain of $\bar{A}$) is a second order linear differential operator, which generates a strongly continuous semi-group $(S_A(t))_{t\geq 0}$ on the space $H^1(\Omega)$ and is self-adjoint with compact resolvent.

◇The considered system is described by the following parabolic partial differential equations

$$\begin{cases} \frac{\partial x}{\partial t}(\xi,t) = Ax(\xi,t) + Bu(t) & Q \\ x(\xi,0) = x_0(\xi) & \Omega \\ \frac{\partial x}{\partial v}(\eta,t) = 0 & \Theta \end{cases} \quad (1)$$

where $\xi \in \Omega, \eta \in \partial\Omega$, $t \in [0,T]$ and $(\xi,t) \in Q, (\eta,t) \in \Theta$.

◇ The measurements can be obtained by using internal zone sensors and pointwise may be located inside $\Omega$ [6-7]. Thus, the output function augmented is given by

$$y(.,t) = Cx(.,t) \quad (2)$$

where the operators $B \in L(R^p, H^1(\Omega))$ and $C \in L(H^1(\Omega), R^q)$ are depended on the structure of actuators and sensors [10-15] (see (figure 2)).

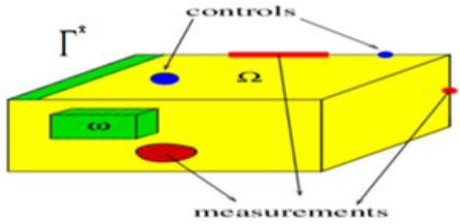

**Fig. 2:** Domain $\Omega$, boundary $\Gamma^*$, and sensors and controls.

◇Under the given assumption above, the system (1) has a unique solution given by the following form [4-5].

$$x(\xi,t) = S_A(t)x_0(\xi) + \int_0^t S_A(t-s)Bu(s)ds \quad (3)$$

◇ The initial state $x_0$ and its gradient $\nabla x_0$ are supposed to be unknown, the problem concerns the reconstruction of the initial gradient $\nabla x_0$ on the region $\Gamma$ of the system domain $\partial\Omega$.

◇ The mathematical modelling illustrated in Figure 2 is more general and complicated than the real modelling in figure 1.

◇The problem is how to construct a an observer to the gradient of current state in a given $\Gamma^*$ may be called $\Gamma^*$AGFO-observer using internal region approach.

◇For deriving $\Gamma^*$AGFO-estimator of $Tx(\xi)$ on $\Gamma^*$, we need to consider the following points:

• Now, we consider the operator $K$ given by the form

$K: H^1(\Omega) \to L^2(0,T,R^q)$
$x \to CS_A(.)x$

where $K$ is bounded linear operator as in [4]. Thus, the adjoint operator $K^*$ of $K$ is defined by $K^*: \mathcal{O} \to X$, and represented by the form

$K^*: L^2(0,T,R^q) \to H^1(\Omega)$
$y^* \to \int_0^t S_A^*(s)C^*y^*(s)ds$

• The operator $\nabla$ denotes the gradient is given by

$\begin{cases} \nabla: H^1(\Omega) \to (H^1(\Omega))^n \\ x \to \nabla_x = (\frac{\partial x}{\partial \xi_1}, \ldots, \frac{\partial x}{\partial \xi_n}) \end{cases}$

with the adjoint of $\nabla$ denotes by $\nabla^*$ is given by [27]

$\begin{cases} \nabla^*: (H^1(\Omega))^n \to H^1(\Omega) \\ x \to \nabla_x^* = v \end{cases}$

where $v$ is a solution of the Dirichlet problem

$\begin{cases} \Delta_v = -div(z) & \Omega \\ v = 0 & \partial\Omega \end{cases}$

• The trace operator of order zero is described by [28]

$\gamma_0: H^1(\Omega) \to H^{1/2}(\partial\Omega)$

which is linear, subjective and continuous [29]. Thus, the extension of the trace operator of order zero which is denoted by $\gamma$ defined as

$\gamma: (H^1(\Omega))^n \to (H^{1/2}(\partial\Omega))^n$

and the adjoints are respectively given by $\gamma_0^*$ and $\gamma^*$.

• For a sub-boundary $\Gamma^*$ of $\partial\Omega$ and let $\tilde{\chi}_{\Gamma^*}$ be the function defined by

$\begin{cases} \tilde{\chi}_{\Gamma^*}: H^{1/2}(\partial\Omega) \to H^{1/2}(\Gamma^*) \\ x \to \chi_{\Gamma^*}x = x|_{\Gamma^*} \end{cases}$

with $x|_{\Gamma^*}$ is the restriction of the state $x$ to $\Gamma^*$, and

$\chi_{\Gamma^*}: (H^1(\partial\Omega))^n \to (H^{1/2}(\Gamma^*))^n$

where the adjoints are respectively given by $\tilde{\chi}_{\Gamma^*}^*$ and $\chi_{\Gamma^*}^*$.

•Finally, we introduced the operator $\chi_{\Gamma^*}\gamma\nabla K^*$ from $\mathcal{O}$ into $(H^{1/2}(\Gamma^*))^n$ and the adjoint of this operator given by $K\nabla^*\gamma^*\chi_{\Gamma^*}^*$.

• We first recall a sensors are defined by any couple $(D,f)$ where $D$ be a non-empty closed subset of $\Omega$ which is represented the spatial supports of sensor and $f \in L^2(D)$ represent the distributions of the sensing measurements on $D$.

Then, according to the choice of the parameters $D$ and $f$, we have different types of sensor:

• It may be zone, if $D \subset \Omega$ and $f \in L^2(D)$. In this case, the operator $C$ is bounded [8-9] and the output function (2) may be given by the form

$$y(t) = \int_D f(\xi)x(\xi,t)d\xi = Cx(\xi,t) \quad (4)$$

• It may be pointwise, if $D = \{b\}$ with $b \in \Omega$ and $f = \delta(.-b)$, where $\delta$ is the Dirac mass concentrated in $b$. In this case, the operator $C$ is un bounded and the output function (2) may be given by the form

$$y(t) = \int_\Omega x(\xi,t)\delta_b(\xi-b)d\xi \quad (5)$$

In this section, we present some definitions and descriptions of regional boundary gradient observability, detectability and strategic sensor, which is derived of [17-22]. Consider the autonomous system of (1) define by

$$\begin{cases} \frac{\partial x}{\partial t}(\xi,t) = Ax(\xi,t) & Q \\ x(\xi,0) = x_0(\xi) & \Omega \\ \frac{\partial x}{\partial v}(\eta,t) = 0 & \Theta \end{cases} \quad (6)$$

The solution of (7) is given by the following form

$$x(\xi,t) = S_A(t)x_0(\xi) \text{ for all } t \in [0,T] \quad (7)$$

• The systems (6)-(7) are said to be exactly regionally boundary gradient observable on $\Gamma^*$ ($E\Gamma^*G$-observable) if

$Im\,\chi_{\Gamma^*}\nabla K^* = (H^{1/2}(\Gamma^*))^n$

• The systems (6)-(7) are said to be weakly regionally boundary gradient observable on $\Gamma^*$ ($W\Gamma^*G$-observable) if

$\overline{Im\chi_{\Gamma^*}\nabla K^*} = (H^{1/2}(\Gamma^*))^n$

It is equivalent to say that the systems (6)-(2) are $W\Gamma^*G$-observable if

$Ker K\nabla^*\chi_{\Gamma^*} = \{0\}$

• If the systems (6)-(2) are is $W\Gamma^*G$-observable, then $x_0(\xi,0)$ is given by

$$x_0 = (K^*K)^{-1}K^*y = K^\dagger y, \quad (8)$$

where $K^\dagger$ is the pseudo-inverse of the operator $K$ [9-10].

• A sensor $(D,f)$ is regional boundary gradient strategic on $\Gamma^*$ ($\Gamma^*G$-strategic) if the observed system is $W\Gamma^*G$-observable.

• The measurements can be obtained by the use of zone or pointwise sensors, which may be located in $\Omega$ [11-12].

## 3. $\Gamma^*G$-obervability and $\Gamma^*AG$-detectability

This section links $\Gamma^*G$-obervability and $\Gamma^*AG$-detectability notions and which roll paly to build the devoted observer.

• The semi-group $(S_A(t))_{t\geq 0}$ is regionally boundary asymptotically gradient stable on $(H^{1/2}(\Gamma^*))^n$ ($\Gamma^*AG$-stable), then for all $x_o \in H^1(\Omega)$, the solution of autonomous system associated to system (1) coverage to zero when $t$ tend to $\infty$.

• The system (6) is said to be $\Gamma^*AG$-stable if the operator $A$ generates a semi-group which is $\Gamma^*AG$-stable.

• A system is said to be $\Gamma^*AG$-stable if and only if there exists some positive constants $M_{\Gamma^*}$, $\alpha_{\Gamma^*}$, such that
$$\|\chi_{\Gamma^*}\gamma\nabla S_A(.)\|_{L\left((H^{1/2}(\Gamma^*))^n, H^1(\Omega)\right)} \leq M_{\Gamma^*}e^{\alpha_{\Gamma^*}}, \forall t \geq 0. \quad (9)$$

• If the semi-group $(S_A(t))_{t\geq 0}$ is $\Gamma^*AG$-stable, then for all $x_o \in H^1(\Omega)$, the solution of autonomous system (6) associated to system (1) satisfies
$$\|\chi_{\Gamma^*}\gamma\nabla x(.,t)\|_{(H^{1/2}(\Gamma^*))^n} =$$
$$\|\chi_{\Gamma^*}\gamma\nabla S_A(t)x_0\|_{(H^{1/2}(\Gamma^*))^n}$$
$$\leq M_{\Gamma^*}e^{\alpha_{\Gamma^*}}\|x_0\|_{(H^{1/2}(\Gamma^*))^n}$$
and then, we have
$$\lim_{t\to\infty}\|\chi_{\Gamma^*}\gamma\nabla x(t)\|_{(H^{1/2}(\Gamma^*))^n} = 0.$$

• The system (1)-(2) is said to be regionally boundary asymptotically gradient detectable on $\Gamma^*$ ($\Gamma^*AG$-detectable), if there exists an operator $H_{\Gamma^*AG}: R^q \to \left(H^{1/2}(\Gamma^*)\right)^n$, such that the operator $(A - H_{\Gamma^*AG}C)$ generates a strongly continuous semi-group $\left(S_{H_{\Gamma^*AG}}(t)\right)_{t\geq 0}$, which is $\Gamma^*AG$-stable.

**Proposition 3.1:** If the system (6)-(2) is $E\Gamma^*G$-observable, then it is $\Gamma^*AG$-detectable. This results gives the following inequality $\exists k_{\Gamma^*} > 0$, such that
$$\|\chi_{\Gamma^*}\gamma\nabla S_A(.)x\|_{(H^{1/2}(\Gamma^*))^n} \leq k_{\Gamma^*}\|CS_A(.)x\|_{L^2(0,\infty,\mathcal{O})}, \quad (10)$$
for all $x \in \left(H^{1/2}(\Gamma^*)\right)^n$.

**Proof:** The proof of this proposition can be concluded from the results on observability by considering the operator $\chi_{\Gamma^*}\gamma\nabla K^*$ in the following forms [29-30]
1. $Imf \subset Img$.
2. There exists $k > 0$, such that
$$\|f^*x^*\|_{E^*} \leq k\|g^*x^*\|_{F^*}, \text{for all } x^* \in G^*$$
From the right hand said of above inequality $k_{\Gamma^*}\|g^*x^*\|_{F^*}$, there exists $M_{\Gamma^*}, \omega_{\Gamma^*} > 0$ with $k_{\Gamma^*} < M_{\Gamma^*}$, such that
$$k_{\Gamma^*}\|g^*x\|_{F^*} \leq M_{\Gamma^*}e^{-\omega_{\Gamma^*}t}\|x^*\|_{F^*}$$
where $E, F$ and $G$ be a reflexive Banach spaces and $f \in L(E,G)$, $g \in L(F,G)$. If we apply this result, considered
$$E = G = (H^{1/2}(\Gamma^*))^n, F = \mathcal{O}, f = Id_{(H^{1/2}(\Gamma^*))^n}$$
and
$$g = S_A^*(.)\chi_{\Gamma}^*\gamma^*\nabla^*C^*$$
where $S_A(.)$ is a strongly continuous semi-group generates by $A$, which is $\Gamma^*AG$-stable on $\Gamma^*$, then it is $\Gamma^*AG$-detectable on $\Gamma^*$. ∎

Thus, the notion of $\Gamma^*AG$-detectability is a weaker property than the $E\Gamma^*G$-observability [30].

**Remark 3.2:** We show that the characterization result that links an $\Gamma^*AG$-detectable and sensors structures. For that purpose, we assume that the operator $A$ has a complete set of eigenfunctions $H^1(\Omega)$[13] denoted $\varphi_{mj}$ orthonormal in $(H^{1/2}(\Gamma^*))^n$ and the associated eigenvalues $\lambda_m$ are of multiplicity $r_m$ and suppose that the system (1) has $J$ unstable modes.

Thus, the sufficient condition of an $\Gamma^*AG$-detectability is given by the following result.

**Theorem 3.3:** Suppose that there are $q$ zone sensors $(D_i, f_i)_{1\leq i\leq q}$ and the spectrum of $A$ contains $J$ eigenvalues with non-negative real parts. The system (1)-(2) are $\Gamma^*AG$-detectable if and only if
1. $q \geq m$,
2. rank $G_i = m_i$, for all $i = 1,...,J$ with
$$G = (G)_{ij} = \begin{cases} \langle\psi_j(.), f_i(.)\rangle_{L^2(D_i)} & \text{zone sensors} \\ \psi_j(b_i) & \text{pointwise sensors} \end{cases}$$
where $\sup m_i = m < \infty$ and $j = 1,...,\infty$.

**Proof:** The proof can be stated by the same way as in ref. [31] with some modifications by choosing pointwise sensors. ∎

3.1. is complete

## 4. Regional internal and $\Gamma^*AGFO$-observer reconstruction

In this section, we give the sufficient conditions which are guarantee the existence of ($\Gamma^*AGFO$-Observer) which allows to construct a $\Gamma^*AGFO$-estimator of the state $\chi_{\Gamma^*}\gamma\nabla Tx(\xi,t)$ by using internal region to pass on the regional boundary. The original results are presented and examined.

### 4.1 Definitions and characterizations
This subsection related to present some definitions and characterizations.

**Definition 3.1:** The dynamical system associated to the considered systems (1)-(2) is given by
$$\begin{cases} z(\xi,t) = F_{\Gamma^*AG}z(\xi,t) + G_{\Gamma^*AG}u(t) + H_{\Gamma^*AG}y(t) & Q \\ z(\xi,0) = z_0(\xi) & \Omega \\ \frac{\partial x}{\partial \nu}(\eta,t) = 0 & \Theta \end{cases} \quad (11)$$
where $F_{\Gamma^*AG}$ generates a strongly continuous semi-group $(S_{F_{\Gamma^*AG}}(t))_{t\geq 0}$ which is $F_{\Gamma^*AG}$-stable on $Z$ and $G_{\Gamma^*AG} \in L(U,Z), H_{\Gamma^*AG} \in L(\mathcal{O},Z)$. The system (11) defines an $\Gamma^*AG$-estimator for $T_{\Gamma^*AG}x(\xi,t) = \chi_{\Gamma}\nabla Tx(\xi,t)$, where $T: X \to Z$ with
$$T_{\Gamma^*AG}x(\xi,t) = z(\xi,t)$$

**Definition 4.2:** Suppose there exists a dynamical system with state $z(.,t) \in Z$ given by
$$\begin{cases} \frac{\partial z}{\partial t}(\xi,t) = Az(\xi,t) + Bu(t) - H_{\Gamma^*AG}C(x(\xi,t) - z(\xi,t)) & Q \\ z(\xi,0) = z_o(\xi) & \Omega \\ z(\eta,t) = 0 & \Sigma \end{cases} \quad (12)$$
In this case the operator $F_{\Gamma^*AG}$ in system (1)[13] is given by $F_{\Gamma^*AG} = A - H_{\Gamma^*AG}C$ where $T_{\Gamma^*AG} = I_{\Gamma^*AGFO}$ the identity operator. Thus the operator $A - H_{\Gamma^*AG}C$ generate a strongly continuous semi-group $(S_{A-H_{\Gamma^*AG}C}(t))_{t\geq 0}$ on separable Hilbert space $Z$ which is $\Gamma^*AG$-stable.

Thus, $\exists M_{A-H_{\Gamma^*AG}C}, \alpha_{A-H_{\Gamma^*AG}C} > 0$ such that
$$\|S_{A-H_{\Gamma^*AG}C}(.)\| \leq M_{A-H_{\Gamma^*AG}C}e^{-\alpha_{A-H_{\Gamma^*AG}C}t}, \forall t \geq 0.$$
Then, let such that solution of (11) similar to (3)
$$z(\xi,t) = S_{A-H_{\Gamma^*AG}C}(t)z(\xi) + \left[\int_0^t S_{A-H_{\Gamma^*AG}C}(t-\tau)Bu(\tau)H_{\Gamma^*AG}y(\tau)\right]d\tau.$$

**Definition 4.3:** The system (12) defines $\Gamma^*AGFO$-estimator such that
$$z(\xi,t) = \chi_{\Gamma^*}\nabla T_{\Gamma^*AGFO}x(\xi,t) = I_{\Gamma^*AGFO}x(\xi,t) \in (H^{1/2}(\Gamma^*))^n$$
where $x(\xi,t)$ is the solution of the systems (1)-(2), if
$$\lim_{t\to\infty}\|z(.,t) - \chi_{\Gamma^*}\nabla T_{\Gamma^*AGFO}x(\xi,t)\|_{(H^{1/2}(\Gamma^*))^n} = 0,$$
and $\chi_{\Gamma^*}\nabla\chi_{\Gamma^*}\nabla I_{\Gamma^*AGFO}$ maps $D(A)$ into $D(A - H_{\Gamma^*AG}C)$ where $z(\xi,t)$ is the solution of system (12).

**Remark 4.4:** The dynamic system (12) specifies $\Gamma^*AGFO$-observer of the systems given by (1)-(2) if the following holds:
1- There exists
$M_{\Gamma^*AGFO} \in L(R,(H^{1/2}(\Gamma^*)^n))$ and $N_{\Gamma^*AGFO} \in L((H^{1/2}(\Gamma^*))^n)$
such that
$M_{\Gamma^*AGFO}C + N_{\Gamma^*AGFO} = I_{\Gamma^*AGFO}$.
2- $A - F_{\Gamma^*AGFO} = H_{\Gamma^*AGFO}C$ and $G_{\Gamma^*AGFO} = B$.
3- The system (11) defines $\Gamma^*AGFO$-estimator for $x(\xi,t)$.
4- The purpose of $\Gamma^*AGFO$-observer is to provide an approximation to the original system state gradient. This approximation is given by
$$\hat{x}(t) = M_{\Gamma^*AGFO}y(t) + N_{\Gamma^*AGFO}z(t).$$

**Definition 4.5:** The systems (1)-(2) are regionally boundary asymptotically gradient full order observable on $\Gamma^*$ ($\Gamma^*AGFO$-observable), if there exists a dynamic system which is $\Gamma^*AGFO$-observer for this system.

**Remark 4.6:** If a system is $\Gamma^*$AGFO–observable then, the corresponding sensor is $\Gamma^*$AGFO- strategic sensor.

• The regional boundary observer in $\Gamma^*$ may be seen as internal regional observer in $\omega_r$ if we consider the following transformations. Let $\Re$ be the continuous linear extension operator [28], $\Re: (H^{1/2}(\partial\Omega))^n \to H^1(\Omega)^n$ such that
$$\chi_\Gamma \gamma \nabla \Re h(\mu,t) = h(\mu,t), \quad \text{for all } h \in \left(H^{1/2}(\Gamma^*)\right)^n \qquad (13)$$

• Let $r > 0$ is an arbitrary and sufficiently small real and let the sets
$$E = \bigcup_{x \in \Gamma} B(x,r) \text{ and } \overline{\omega}_r = E \cap \overline{\Omega}$$
where $B(x,r)$ is the ball of radius $r$ centered in $x(\xi,t)$ and where $\Gamma^*$ is a part of $\overline{\omega}_r$ (Fig. 3).

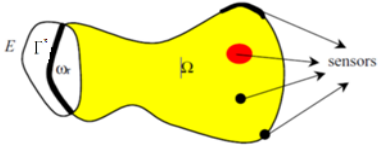

**Fig. 3:** Domain $\Omega$, region $\overline{\omega}_r$ and the sub-region $\Gamma^*$.

• For the sub-region $\omega_r$ of the domain $\Omega$ and let $\chi_{\omega_r}$ be the function defined by
$$\chi_{\omega_r}: (L^2(\Omega))^n \to (L^2(\omega_r))^n$$
$$x \to \chi_{\omega_r} x = x|_{\omega_r}$$
where $x|_{\omega_r}$ is the restriction of $x$ to $\omega_r$ for more derails see ref. s ([10-13]).

**Proposition 4.7:** If the system (1)-(2) is $E\overline{\omega}_r G$-observable (respectively $W\overline{\omega}_r G$), then it is $E\Gamma^* G$ -observable (respectively $W\Gamma^* G$-observable) (see [30]).

From proposition 4.7, we can deduce the following important result.

**Proposition 4.8:** A dynamical system is $\overline{\omega}_r$ AGFO-observer for the system (1)-(2), then it is $\Gamma^*$AGFO-observer.

**Proof:** Let $z(\xi,t) \in (H^{1/2}(\Gamma))^n$, $\chi^*_{\Gamma^*} z(\xi,t) \in (H^{1/2}(\Omega))^n$ and $\overline{z}(\xi,t)$ be an extension to $(H^{1/2}(\partial\Omega))^n$. By using equation (10) and trace theorem there exists
$$\Re \chi^*_{\Gamma^*} z(\xi,t) \in (H^1(\Omega))^n,$$
with bounded support [15] such that
$$\gamma \nabla (\Re \chi^*_{\Gamma^*} z(\xi,t)) = \overline{z}(\xi,t) \qquad (14)$$
Since the system (12) is regional $\overline{\omega}_r AGFO$-observer, then it is $\omega_r AGFO$-observer, there exists a dynamical system with $x(\xi,t) \in X$ such that
$$\chi_{\omega_r} \nabla Tx(\xi,t) = \chi_{\omega_r} \Re \chi^*_{\Gamma^*} z(\xi,t)$$
then we have
$$\chi_\Gamma (\gamma \chi^*_\omega \chi_\omega \nabla Tx)(\xi,t) = x(\xi,t) \qquad (15)$$
The equations (2) and (15) allow
$$\begin{bmatrix} y(\xi,t) \\ z(\xi,t) \end{bmatrix} = \begin{bmatrix} C \\ \chi_\Gamma(\chi^*_\omega \chi_\omega \nabla T) \end{bmatrix} x(\xi,t)$$
and there exist two linear bounded operators $\overline{R}$ and $\overline{S}$ satisfy the relation
$$\overline{R}C + \chi_{\Gamma^*}(\gamma \chi^*_\omega \chi_\omega \nabla T) = I_{\Gamma^*}$$
There exists an operator $F_{\overline{\omega}r}$ is regionally stable on $\overline{\omega}_r$, then it is regionally stable on $\Gamma^*$ [21-20].
Finally the system (12) is a $\Gamma^*$AGFO-observer [30-31]. ∎

## 4.2 Sufficient condition for $\Gamma^*$AGFO-observer

As in (Refs. [21, 30-31]), we extend the characterization result that links the $\Gamma^*$AGFO-observer and $\Gamma^*$AGFO-detectability which is described a sufficient condition for $\Gamma^*$AGFO-observer in the following main result.

**Theorem 4.8:** If the system (1)-(2) is $\Gamma^*$AGFO-detectable, then, the dynamical system (12) is the associated $\Gamma^*$AGFO-observer, i.e.
$$\|x(.,t) - T_{\Gamma^*AGFO}z(.,t)\|_{H^{\frac{1}{2}}(\Gamma^*)} =$$
$$\lim_{t\to\infty} \|x(.,t) - I_{\Gamma^*AGFO}z(.,t)\|_{H^{1/2}(\Gamma^*)} =$$
$$\lim_{t\to\infty} \|x(.,t) - I_{\Gamma^*AGFO}z(.,t)\|_{H^{1/2}(\Gamma^*)} =$$
$$\lim_{t\to\infty} \|x(.,t) - z(.,t)\|_{H^{1/2}(\Gamma^*)} = 0 \qquad (16)$$

**Proof:** From the assumptions of section 2, the system (1) can be decomposed by the projections $P$ and $I - P$ on two parts, unstable and stable [5]. The state vector may be given by where $x_1(\xi,t)$ is the state component of the unstable part of the system (1), may be written in the form
$$\begin{cases} \frac{\partial x_1}{\partial t}(\xi,t) = Ax_1(\xi,t) + Bu(t) & Q \\ x_1(\xi,t) = x_{0_1}(\xi) & \Omega \\ \frac{\partial x_1}{\partial v}(\eta,t) = 0 & \Theta \end{cases} \qquad (17)$$
and $x_2(\xi,t)$ is the component state of the stable part of the system (1) given b
$$\begin{cases} \frac{\partial x_2}{\partial t}(\xi,t) = Ax_{21}(\xi,t) + Bu(t) & Q \\ x_2(\xi,t) = x_{0_2}(\xi) & \Omega \\ \frac{\partial x_2}{\partial v}(\eta,t) = 0 & \Theta \end{cases} \qquad (18)$$
The operator $A_1$ is represented by a matrix of order $(\sum_{n=1}^J s_n, \sum_{n=1}^J s_n)$ given
$$A_1 = diag[\lambda_1, \ldots, \lambda_1, \ldots, \lambda_2, \ldots, \lambda_2, \ldots, \lambda_J, \ldots, \lambda_J]$$
and
$$PB = [G_1^{tr}, G_2^{tr}, \ldots, G_j^{tr}]$$
Put $e(\xi,t) = x(\xi,t) - z(\xi,t)$ where $z(\xi,t)$ is the solution of the system (12). By deriving the above equation and substituting equations (1) and (12), we obtain
$$\frac{\partial e}{\partial t}(\xi,t) = \frac{\partial x}{\partial t}(\xi,t) - \frac{\partial z}{\partial t}(\xi,t)$$
$$= Ax(\xi,t) - Az(\xi,t) - H_{\Gamma^*AGFO}C(x(.,t) - z(\xi,t))$$
$$= (A - H_{\Gamma^*AGFO}C)e(\xi,t)$$
Since the system (1)-(2) is $\Gamma^*$AG-detectable, there exists an operator $H_{\Gamma^*AGFO} \in L(R^q, H^{1/2}(\Gamma^*))$, such that the operator $(A - H_{\Gamma^*AGFO}C)$, generates a stable, strongly continuous semi-group $(S_{H_{\Gamma^*AGFO}}(t))_{t\geq 0}$ on the space $H^{1/2}(\Gamma^*)$, that means $\exists M_{\Gamma^*AGFO}, \alpha_{\Gamma^*AGFO} > 0$, which is satisfied the following inequality
$$\|\chi_{\Gamma^*}\gamma\nabla S_A(.)\|_{H^{1/2}(\Gamma^*)} \leq M_{\Gamma^*AGFO}e^{-\alpha_{\Gamma^*AGFO}t}$$
Finally, we have
$$\|e(.,t)\|_{H^{1/2}(\Gamma^*)} \leq$$
$$\|\chi_{\Gamma^*}\gamma\nabla S_{H_{\Gamma^*AGFO}}(.)\|_{H^{1/2}(\Gamma^*)}\|e_0(.)\|_{H^{1/2}(\Gamma^*)}$$
$$\leq M_{\Gamma^*AGFO}e^{-\alpha_{\Gamma^*AGFO}t}\|e_0(.)\|_{H^{1/2}(\Gamma^*)}$$
and
$$e_0(\xi,t) = x(\xi,t) - z(\xi,t)$$
therefore
$$\lim_{t\to\infty} \|e(.,t)\|_{H^{1/2}(\Gamma^*)} = 0.$$
Consequently, the dynamical system (12) is a $\Gamma^*$AGFO-observer for the considered system (1)-(2). ∎

**Remark 4.9.** From theorem 4.8., we can deduce the following results:
1. A dynamical system which is an $\partial\Omega$AGFO-observer is $\Gamma^*$AGFO-observer.
2. If a system is $\Gamma_1^*$AGFO-observer, then it is $\Gamma_2^*$AGFO-observer in every $\Gamma_1^* \subset \Gamma_2^*$, but the converse is not true. This may be proven in the following application

## 4.3 Counter Application to Diffusion System

Consider two-dimensional of system (1) which is given by the following diffusion parabolic equation*s*
$$\begin{cases} \frac{\partial x}{\partial t}(\xi_1,\xi_2,t) = \Delta x(\xi_1,\xi_2,t) + \delta_{\overline{b}}(\xi_1,\xi_2)u(t) & Q \\ x(\xi_1,\xi_2,t) = x_0(\xi_1,\xi_2) & \Omega \\ \frac{\partial x}{\partial v}(\eta_1,\eta_1,t) = 0 & \Theta \\ y(t) = \delta_b(b_1,b_2)x(\xi_1,\xi_2,t) & Q \end{cases} \qquad (19)$$
where
$$\Omega = ]0,1[\times]0,1[, \delta_{\overline{b}}(\xi_1,\xi_2) = \delta(\xi_1 - \overline{b}_1, \xi_2 - \overline{b}_1)$$

and $\bar{b} = (b_{\bar{b}_1}, \bar{b}_2) \in \Omega$ is location of the internal pointwise control $(\bar{b}, \delta_{\bar{b}})$. Then, and the operator $Bu(t)$ in system (19) is given by
$$Bu(t) = \delta_{\bar{b}}(\bar{b}_1, \bar{b}_2)u(t) \quad (20)$$
Consider the internal filament sensor where σ = $Im(\gamma) \subset \Omega$ is symmetric with respect to the line $b = (b_1, b_2)$ as in and $f(b_1, b_2) = \cos \pi b_1 \cos \pi b_2$.
The augmented output function (2) can be written by
$$y(t) = \int_\Omega \delta_b(\xi_1 - b_1, b_2 - \xi_2) x(\xi_1, \xi_2, t) d\xi_1 d\xi_2 \quad (21)$$
Since the pointwise sensor is a couple $(b, \delta_b)$ of $b$ and $\delta_b$, then
$$y(t) = Cx(\xi_1, \xi_2, t) = \delta_b(b_1, b_2)x(\xi_1, \xi_2, t) \quad (22)$$
The operator $A = \Delta$ generates a strongly continuous semi-group $(S_A(t))_{t \geq 0}$ on the Hilbert space $H^1(\Omega)$ given by
$$S_A(t)x = \sum_{n,m=0}^{\infty} e^{\lambda_{nm} t} \langle x, \varphi_{nm} \rangle_{H^1(\Omega)} \varphi_{nm} \quad (23)$$
where
$$\lambda_{nm} = -(n^2, m^2)\pi^2 \quad (24)$$
are the eigenvalues and
$$\varphi_{nm}(\xi_1, \xi_2) = 2a_{nm} \cos(n\pi\xi_1) \cos(n\pi\xi_2) \quad (25)$$
are eigenvectors with
$$2a_{nm} = (1 - \lambda_{nm})^{-1/2} \quad (26)$$
Consider now, the dynamical system
$$\begin{cases} \frac{\partial z}{\partial t}(\xi_1, \xi_2, t) = \Delta z(\xi_1, \xi_2, t) + \delta_{\bar{b}}(\xi_1, \xi_2)u(t) \\ -H\delta_b(b_1, b_2)(x(\xi_1, \xi_2, t) - z(\xi_1, \xi_2, t)) & Q \\ z(\xi_1, \xi_2, 0) = z_0(\xi_1, \xi_2) & \Omega \\ \frac{\partial z}{\partial \nu}(\eta_1, \eta_1, t) = 0 & \Theta \end{cases} \quad (27)$$
where $H \in L(R^q, Z), Z$ is a Hilbert space and $C: H^1(\bar{\Omega}) \to R^q$ is a linear operator. If the state $x_0$ is defined in Ω by
$$x_0(\xi_1, \xi_2) = \cos(\pi\xi_1) \cos(2\pi\xi_2), \quad (28)$$
then the system (19) is not WΩG-observable, i.e. $(\sigma, f)$ is not ΩG-strategic sensor [17-18]and therefore the system (19) is not ΩAG-detectable.

Thus, the dynamical system (27) is not ΩAGFO-observer [30] for the system (19 (see [21]). Now, consider the regionΓ* = ]0,1[× {1} ⊂ ∂Ω (figure 4) with previous results, then

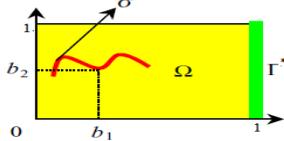

**Fig. 4:** Domain Ω, region Γ* and locations σ of filament pointwise sensor.

The dynamical system
$$\begin{cases} \frac{\partial z}{\partial t}(\xi_1, \xi_2, t) = \Delta z(\xi_1, \xi_2, t) + \delta_{\bar{b}}(\xi_1, \xi_2)u(t) \\ -H_{\Gamma^*AGFO}\delta_b(b_1, b_2)(z(\xi_1, \xi_2, t) - x(\xi_1, \xi_2, t)) & Q \\ z(\xi_1, \xi_2, 0) = z_0(\xi_1, \xi_2) & \Omega \\ \frac{\partial z}{\partial \nu}(\eta_1, \eta_1, t) = 0 & \Theta \end{cases} \quad (29)$$
where $H \in L(R^q, H^{1/2}(\Gamma^*))$, then in this case, the system (19) is WΓ*G-observable. Thus, the sensor $(\sigma, \delta_b)$ is Γ*G-strategic [20] if,
$nb_1 \notin N$ and $nb_2 \notin N$ for every $n, m = \{1, ..., J\}$.
Hence, the system (19) is Γ*AG-detectable [1]. Finally, the dynamical system (27) isΓ*AGFO-observer for the system (19) [21].∎

**Remark 4.11:** If the system is Γ*AG-detectable, then it is possible to construct Γ*AGFO-observer for the original system.

**Acknowledgements.** Our thanks in advance to the editors and experts for considering this paper to publish in this esteemed journal. The authors appreciate your time and effort in reviewing the manuscript and greatly value your assistant as reviewer for the paper.

## 6. Conclusion

We have extended the original results related to the concept of Γ*AGFO-observerfor parabolic distributed system in where the dynamic system generates a strongly continuous semi-group Hilbert space.
More precisely, we have shown that, the possibility to design a dynamic system which is enable to observe asymptotically the state gradient in sub-region Γ* of the boundary ∂Ω using the corresponding detectability and strategic sensors in different situations.
The problem of passage form internal region to regional boundary case is proved and analyzed with an application to diffusion system. Moreover, many problem still opened like the development of these results to case of hyperbolic distributed parameter systems as in [25].